\documentclass[12pt,oneside,reqno]{amsart}
\usepackage{mathrsfs}
\usepackage{graphics}
\usepackage{color, amssymb}
\newcommand{\dif}{\mathrm{d}}
\newcommand{\me}{\mathrm{e}}
\newcommand{\be}{\begin{eqnarray}}
\newcommand{\ee}{\end{eqnarray}}
\newcommand{\ce}{\begin{eqnarray*}}
\newcommand{\de}{\end{eqnarray*}}
\newtheorem{theorem}{Theorem}[section]
\newtheorem{lemma}[theorem]{Lemma}
\newtheorem{remark}[theorem]{Remark}
\newtheorem{definition}[theorem]{Definition}
\newtheorem{proposition}[theorem]{Proposition}
\newtheorem{example}[theorem]{Example}
\newtheorem{corollary}[theorem]{Corollary}

\newcommand{\cN}{{\mathcal{N}}}

\newcommand{\cL}{{\mathcal{L}}}

\newcommand{\PX}{{\Bbb{P}}}

\def\e{\varepsilon}

\def\a{\alpha}

\def\p{\partial}

\def\[{{\Big[}}
\def\]{{\Big]}}
\def\<{{\langle}}
\def\>{{\rangle}}
\def\({{\Big(}}
\def\){{\Big)}}

\def\no{\nonumber}
\def\bt{\begin{theorem}}
\def\et{\end{theorem}}
\def\bl{\begin{lemma}}
\def\el{\end{lemma}}
\def\br{\begin{remark}}
\def\er{\end{remark}}
\def\bx{\begin{Examples}}
\def\ex{\end{Examples}}
\def\bd{\begin{definition}}
\def\ed{\end{definition}}
\def\bp{\begin{proposition}}
\def\ep{\end{proposition}}
\def\bc{\begin{corollary}}
\def\ec{\end{corollary}}

\def\cB{{\mathcal B}}
\def\cC{{\mathcal C}}
\def\cD{{\mathcal D}}

\def\cF{{\mathcal F}}
\def\cG{{\mathcal G}}

\def\cI{{\mathcal I}}

\def\cL{{\mathcal L}}

\def\cN{{\mathcal N}}

\def\mE{{\mathbb E}}

\def\mP{{\mathbb P}}

\def\mR{{\mathbb R}}

\def\geq{\geqslant}
\def\leq{\leqslant}

\begin{document}

\allowdisplaybreaks

\title{Escape probability for stochastic dynamical systems with jumps*}

\author{Huijie Qiao$^1$, Xingye Kan$^2$  and Jinqiao Duan$^2$}



\thanks{*This work is partially supported by the NSF of China (No. 11001051 and No. 11028102) and the NSF grant DMS-1025422.}


\date{\today (Revised version)}

\dedicatory{1. Department of Mathematics,
Southeast University\\
Nanjing, Jiangsu 211189,  China\\
hjqiaogean@yahoo.com.cn \\
2. Institute for Pure and Applied Mathematics, University of California\\
Los Angeles, CA 90095, USA\\ \& \\Department of Applied Mathematics, Illinois Institute of Technology\\
Chicago, IL 60616, USA\\
xkan@hawk.iit.edu,  \;\;  duan@iit.edu }

\begin{abstract}
The escape probability is a deterministic concept that quantifies some aspects of stochastic dynamics.
This issue has been investigated previously for dynamical systems driven by   Gaussian Brownian motions.  The present work considers escape probabilities for  dynamical systems
driven by non-Gaussian L\'evy motions, especially symmetric $\alpha$-stable L\'evy motions. The escape probabilities are characterized as solutions of the
   Balayage-Dirichlet problems of certain partial differential-integral equations. Differences between escape probabilities for dynamical systems driven by Gaussian and non-Gaussian noises are highlighted.  In certain special cases, analytic results for escape probabilities are given.
\end{abstract}

\maketitle \rm

{\it Dedicated to Professor David Nualart on the occasion of his 60th birthday}

\section{Introduction}

Stochastic dynamical systems arise as mathematical models for complex phenomena in biological, geophysical, physical and chemical sciences, under random fluctuations. A specific orbit (or trajectory) for such a system could vary wildly from one realization to another, unlike the situation for deterministic dynamical systems. It is desirable to have different concepts for quantifying stochastic dynamical behaviors.
The escape probability is such a concept.

Brownian motions are  Gaussian stochastic processes and thus are appropriate for modeling Gaussian random fluctuations.
Almost all sample paths of Brownian motions are continuous in time.
For a dynamical system  driven by Brownian motions, almost all orbits (or paths or trajectories) are thus continuous in time. The \emph{escape probability}   is the likelihood that
  an orbit,  starting inside an open domain $D$,   exits this domain first through a specific part $\Gamma $ of the boundary $\p D$. This concept helps understand   various phenomena in sciences. One example is in molecular genetics   \cite{smith}. The frequency of collisions of two
single strands of long helical DNA molecules that leads to a
double-stranded molecule is of interest and can be computed by virtu
of solving an escape probability problem. It turns out that the escape probability satisfies an elliptic
partial differential equation with properly chosen boundary conditions \cite{smith,clls, Schuss, BrannanDuanErvin}.


Non-Gaussian random fluctuations are widely observed in
various areas such as physics, biology, seismology, electrical
engineering and finance  \cite{Woy, Koren, Mat}.   L\'evy motions are a large class of non-Gaussian stochastic processes whose sample paths are discontinuous in time.  For a dynamical system  driven by L\'evy motions,  almost all the orbits $X_t$ are discontinuous in time. In fact, these orbits are c\`adl\`ag (right-continuous with left limit at each time instant), i.e., each of these orbits has countable jumps in time. Due to these jumps, an orbit could escape an open domain   without passing through its boundary. In this case, the \emph{escape probability} is the likelihood that an
 orbit,  starting inside an open domain $D$,   exits this domain first by landing in a target domain $U$ in $D^c$ (the complement of domain $D$).

  As we see, the escape probability is defined slightly differently for dynamical systems driven by Gaussian or non-Gaussian processes. Although the escape probability for the former has been investigated extensively,
 the characterization for the escape probability  for the latter has not been well documented as a dynamical systems analysis tool for applied  mathematics and   science  communities. See our recent works \cite{ChenDuan, GaoDuan} for numerical analysis of escape probability and mean exit time for dynamical systems driven by symmetric $\alpha$-stable L\'evy motions.


In this paper, we carefully derive a partial differential-integral equation to be satisfied by the escape probability for a class of dynamical systems driven by L\'evy motions, especially symmetric $\alpha$-stable L\'evy motions. Namely the escape probability is a solution of a nonlocal differential equation.
We highlight the differences between escape probabilities for dynamical systems driven by Gaussian and non-Gaussian processes. These are illustrated in a few   examples.


\bigskip

More precisely, let $\{X_t, t\geq0\}$ be a $\mR^d$-valued Markov
process defined on a complete filtered probability space
$(\Omega,\cF,\{\cF_t\}_{t\geq 0}, \PX)$. Let $D$ be an open domain
in $\mR^d$. Define the \emph{exit time}
 \ce
 \tau_{D^c}:=\inf\{t>0: X_t\in D^c\},
 \de
where $D^c$ is the complement of $D$ in $\mR^d$. Namely, $\tau_{D^c}$ is the  first time when $X_t$ hits $D^c$.

\medskip

When $X_t$ has almost surely continuous paths, i.e., $X_t$ is either a Brownian motion or a solution process for a dynamical system  driven by   Brownian motions,   a path starting
at $x\in D$ will hit $D^c$ by hitting $\partial D$ first (assume for
the moment that $\partial D$ is smooth). Thus
$\tau_{D^c}=\tau_{\partial D}$. Let $\Gamma$ be a subset of the boundary $\partial
D$. The likelihood  that $X_t$, starting at $x$, exits from $D$ first through $\Gamma$
is called the escape probability from $D$ to $\Gamma$, denoted as
$p(x)$. That is,
\ce p(x)=\PX \{X_{\tau_{\partial D}}\in\Gamma\}.
\de
We will verify that (Section \ref{sde-bm}) the escape probability $p(x)$ solves
the following Dirichlet boundary value problem:
\be
\left\{\begin{array}{l}
\cL p=0, \;\; x\in D, \\
p|_{\partial D}=\psi,
\end{array}
\right.
\label{dirpro1}
\ee
where $\cL$ is the infinitesimal generator of the process $X_t$ and the boundary data
  $\psi$ is defined as follows
\ce
\psi(x)=\left\{\begin{array}{l}
1, \quad x\in \Gamma,\\
0, \quad x\in \partial D\setminus\Gamma.
\end{array}
\right.
\de

\bigskip

When $X_t$ has c\`adl\`ag
paths which have countable jumps in time, i.e., $X_t$ could be either a L\'evy motion or a solution process of a dynamical system driven by L\'evy motions,   the first hitting of $D^c$ may occur
somewhere in $D^c$. For this reason, we take a
subset $U$ of $D^c$, and define the likelihood that   $X_t$ exits
firstly from $D$ by landing in the target set $U$ as the escape probability from $D$ to
$U$, also denoted by $p(x)$. That is,
\ce
p(x)=\PX \{X_{\tau_{D^c}}\in U\}.
\de

We will demonstrate that (Section \ref{sde-levy}) the escape probability $p(x)$
solves the following Balayage-Dirichlet boundary value problem:
\be
\left\{\begin{array}{l}
Ap=0,\;\; x\in D,\\
p|_{D^c}=\varphi,
\end{array}
\right.
\label{dirpro2}
\ee
where $A$ is the characteristic operator of $X_t$ and $\varphi$ is defined as follows
\ce
\varphi(x)=\left\{\begin{array}{l}
1, \quad x\in U,\\
0, \quad x\in D^c\setminus U.
\end{array}
\right.
\de

Therefore by solving a deterministic boundary value problem
\eqref{dirpro1} or \eqref{dirpro2}, we obtain the escape probability
$p(x)$.

\medskip

This paper is arranged as follows. In Section \ref{prelim}, we introduce Balayage-Dirichlet
problem for discontinuous Markov processes, and also define L\'evy motions. The main result
is stated and proved in Section \ref{BVP}. In Section \ref{example}, we present analytic solutions    for escape probabilities   in  a few special cases.

\section{Preliminaries}\label{prelim}

In this section, we recall  basic concepts and results that will
be needed throughout the paper.

\subsection{Balayage-Dirichlet problem for discontinuous Markov processes}

The following materials are from \cite{bg, l, song93, c, mzz, gm}.
Let $\cG$ be a locally compact space with a countable base and $\mathscr{G}$
be the Borel $\sigma$-field of $\cG$. Also, $\varsigma$ is adjoined to $\cG$ as the
point at infinity if $\cG$ is noncompact,
and as an isolated point if $\cG$ is compact. Furthermore, let $\mathscr{G}_{\varsigma}$
be the $\sigma$-field of Borel sets of $\cG_{\varsigma}=\cG\cup\{\varsigma\}$.

\bd\label{hunt process}
A Markov process $Y$ with state space $(\cG,\mathscr{G})$ is called a Hunt process provided:

(i) The paths functions $t\rightarrow Y_t$ are right continuous on $[0,\infty)$
and have left-hand limits on $[0,\zeta)$ almost surely, where $\zeta:=
\inf\{t:Y_t=\varsigma\}$.

(ii) $Y$ is strong Markov.

(iii) $Y$ is quasi-left-continuous: whenever $\{\tau_n\}$ is an increasing sequence
of $\cF_t$-stopping times with limit $\tau$, then almost surely $Y_{\tau_n}\rightarrow Y_{\tau}$
on $\{\tau<\infty\}$.
\ed

\begin{definition}
Let $G$ be an open subset of $\cG$ and $Y_t(x)$ be a Hunt process starting at $x\in G$.
A nonnegative function $h$ defined on $\cG$ is said to be harmonic with respect to $Y_t$ in $G$ if for every compact set $K\subset G$,
\ce
\mE[h(Y_{\tau_{K^c}}(x))]=h(x), \quad x\in G.
\de
\end{definition}

\bd\label{dirpro}
Let $f$ be nonnegative on $G^c$. We say $h$ defined on $\cG$ solves the Balayage-Dirichlet problem
for $G$ with ``boundary value" $f$, denoted by $(G, f)$, if $h=f$ on $G^c$, $h$ is harmonic with respect to $Y_t$ in $G$
and further satisfies the following boundary condition:
\ce
\forall z\in\partial G, \quad h(y)\rightarrow f(z), \quad ~\mbox{as}~y\rightarrow z~\mbox{from inside}~G.
\de
\ed

A point $z\in\partial G$ is called regular for $G^c$ with respect to $Y_t(z)$ if
$$
\mP\{\tau_{G^c}=0\}=1.
$$
Here $G$ is said to be regular if any $z\in\partial G$ is regular for $G^c$.

Let $\rho$ be a metric on $\cG$ compatible with the given topology. Let $\cI_G$
be the family of functions $g\geq0$ bounded on $\cG$ and lower semicontinuous in $G$ such that
$\forall x\in G$, there is a number $Ag(x)$ satisfying
\ce
\frac{\mE[g(Y_{\tau_{\e}}(x))]-g(x)}{\mE[\tau_{\e}]}\rightarrow Ag(x), ~\mbox{as}~\e\downarrow0,
\de
where $\tau_{\e}:=\inf\{t>0: \rho(Y_t(x),x)>\e\}$. We call $A$ with domain $\cI_G$ the
characteristic operator of $Y_t$ relative to $G$. If $\cL$ with domain $\cD_G$ is the
infinitesimal generator of $Y_t$ relative to $G$, $\cD_G\subseteq\cI_G$ and
\ce
Af=\cL f, \quad f\in \cD_G.
\de
(c.f. \cite{dy})

We quote the following result about the existence and regularity of the solution for the Balayage-Dirichlet problem.
\bt \label{soldir} (\cite{l})\\
Suppose that $G$ is relatively compact and regular,  and  $f$ is nonnegative and bounded
on $G^c$. If $f$ is continuous at any $z\in\partial G$, then $h(x)=\mE[f(Y_{\tau_{G^c}}(x))]$ is
the unique solution to the Balayage-Dirichlet problem $(G, f)$, and $Ah(x)=0$ for $h\in I_G$.
\et

\subsection{L\'evy motions}

\bd\label{levy}
A process $L_t $, with $L_0=0$ a.s. is a
$d$-dimensional L\'evy process or L\'evy motion if

(i) $L_t$ has independent increments; that is, $L_t-L_s$ is
independent of $L_v-L_u$ if $(u,v)\cap(s,t)=\emptyset$;

(ii) $L_t$ has stationary increments; that is, $L_t-L_s$ has the same
distribution as $L_v-L_u$ if $t-s=v-u>0$;

(iii) $L_t$ is stochastically continuous;

(iv) $L_t$ is right continuous with left limit.
\ed
The characteristic function for $L_t$ is given by
\ce
\mE\left(\exp\{i\<z,L_t\>\}\right)=\exp\{t\Psi(z)\}, \quad
z\in\mR^d,
\de
where $\<\cdot,\cdot\>$ is the scalar product in $\mR^d$. The function $\Psi: \mR^d\rightarrow\mathcal {C}$ is
called the characteristic exponent of the L\'evy process $L_t$. By the
L\'evy-Khintchine formula, there exist a nonnegative-definite
$d\times d$ matrix $Q$, a measure $\nu$ on $\mR^d$ satisfying \ce
\nu(\{0\})=0 ~\mbox{and}~ \int_{\mR^d\setminus\{0\}}(|u|^2\wedge1)\nu(\dif
u)<\infty, \de and $\gamma\in\mR^d$ such that
\be
\Psi(z)&=&i\<z,\gamma\>-\frac{1}{2}\<z,Qz\>\no\\
&&+\int_{\mR^d\setminus\{0\}}\big(e^{i\<z,u\>}-1-i\<z,u\>1_{|u|\leq1}\big)\nu(\dif
u).
\label{lkf}
\ee
The measure $\nu$ is called the L\'evy measure of $L_t$, $Q$ is the diffusion matrix, and $\gamma$ is the drift vector.

\medskip

We now introduce a special class of L\'evy motions, i.e., the symmetric $\alpha$-stable L\'evy motions   $L_t^\alpha$.

\bd\label{rid}
For $\alpha\in(0,2)$. A $d$-dimensional symmetric $\alpha$-stable L\'evy motion
  $L_t^\alpha$ is a L\'evy process with characteristic exponent
\be
\Psi(z)=-C|z|^\alpha, \quad z\in\mR^d,
\label{stablece}
\ee
where
$$
C=\pi^{-1/2}\frac{\Gamma((1+\a)/2)\Gamma(d/2)}{\Gamma((d+\a)/2)}.
$$
(c.f. \cite[Page 115]{s} for the above formula of $C$.)
\ed

Thus, for a $d$-dimensional symmetric $\alpha$-stable L\'evy motion  $L_t^\alpha$, the diffusion matrix $Q=0$, the drift vector $\gamma=0$,  and the
L\'evy measure $\nu$ is given by
$$
\nu(\dif u)=\frac{C_{d,\alpha}}{|u|^{d+\alpha}}\dif u,
$$
where
$$
C_{d,\alpha}=\frac{\a\Gamma((d+\a)/2)}{2^{1-\a}\pi^{d/2}\Gamma(1-\a/2)}.
$$
(c.f. \cite[Page 1312]{cks} for the above formula of $C_{d,\alpha}$.) Moreover,
comparing (\ref{stablece}) with (\ref{lkf}), we obtain
\ce
-C|z|^\alpha=\int_{\mR^d\setminus\{0\}}\big(e^{i\<z,u\>}-1-i\<z,u\>1_{|u|\leq1}\big)
\frac{C_{d,\alpha}}{|u|^{d+\alpha}}\dif u.
\de

Let $\cC_0(\mR^{d})$ be the space of continuous functions $f$ on $\mR^{d}$ satisfying
$\lim\limits_{|x|\rightarrow\infty}f(x)=0$ with norm $\|f\|_{\cC_0(\mR^{d})}=\sup\limits_{x\in\mR^{d}}|f(x)|$.
Let $\cC^2_0(\mR^{d})$ be the set of $f\in\cC_0(\mR^{d})$ such that $f$ is two
times differentiable and the first and second order partial derivatives of $f$ belong
to $\cC_0(\mR^{d})$. Let $\cC_c^{\infty}(\mR^{d})$ stand for the space of all infinitely
differentiable functions on $\mR^{d}$ with compact supports. Define
\ce
(\cL_\alpha f)(x):=\int_{\mR^d\setminus\{0\}}\big(f(x+u)-f(x)
-\<\partial_xf(x),u\>1_{|u|\leq1}\big)\frac{C_{d,\alpha}}{|u|^{d+\alpha}}\dif u
\de
on $\cC_0^2(\mR^{d})$. And then for $\xi\in\mR^{d}$
\ce
(\cL_\alpha e^{i\<\cdot,\xi\>})(x)=e^{i\<x,\xi\>}
\int_{\mR^d\setminus\{0\}}\big(e^{i\<u,\xi\>}-1-i\<\xi,u\>1_{|u|\leq1}\big)\frac{C_{d,\alpha}}{|u|^{d+\alpha}}\dif u.
\de
By Courr\`ege's second theorem (\cite[Theorem 3.5.5, p.183]{da}),
for every $f\in\cC_c^\infty(\mR^{d})$
\ce
(\cL_\alpha f)(x)&=&\frac{1}{(2\pi)^{d/2}}\int_{\mR^d}e^{i\<z,x\>}\left[e^{-i\<x,z\>}
(\cL_\alpha e^{i\<\cdot,z\>})(x)\right]\hat{f}(z)\dif z\\
&=&\frac{1}{(2\pi)^{d/2}}\int_{\mR^d}e^{i\<z,x\>}\left[\int_{\mR^d\setminus\{0\}}
\big(e^{i\<u,z\>}-1-i\<z,u\>1_{|u|\leq1}\big)\frac{C_{d,\alpha}}{|u|^{d+\alpha}}\dif u\right]\hat{f}(z)\dif z\\
&=&-\frac{C}{(2\pi)^{d/2}}\int_{\mR^d}e^{i\<z,x\>}|z|^\alpha\hat{f}(z)\dif z\\
&=&C\cdot[-(-\Delta)^{\alpha/2}f](x).
\de

Set $p_t:=L_t-L_{t-}$. Then $p_t$ defines a stationary $(\cF_t)$-adapted Poisson point process with values in
$\mR^d\setminus\{0\}$ (\cite{iw}). And the characteristic measure of $p$ is the L\'evy measure $\nu$.
Let $N_{p}((0,t],\dif u)$ be the counting measure of
$p_{t}$, i.e., for $B\in\cB(\mR^d\setminus\{0\})$
$$
N_p((0,t],B):=\#\{0<s\leq t: p_s\in B\},
$$
where $\#$ denotes the cardinality of a set. The compensator measure
of $N_p$ is given by
$$
\tilde{N}_{p}((0,t],\dif u):=N_{p}((0,t],\dif u)-t\nu(\dif u).
$$
The L\'evy-It\^o theorem states that for a symmetric $\alpha$-stable process $L_t$,

(i) for $1\leq\alpha<2$,
\ce
L_t=\int_0^t\int_{|u|\leq1}u\tilde{N}_{p}(\dif s, \dif u)
+\int_0^t\int_{|u|>1}uN_{p}(\dif s, \dif u),
\de

(ii) for $0<\alpha<1$,
\ce
L_t=\int_0^t\int_{\mR^d\setminus\{0\}}uN_{p}(\dif s, \dif u).
\de

\section{Boundary value problems for escape probability}
\label{BVP}

In this section, we formulate boundary value problems for the escape probability
associated with Brownian motions, \textsc{SDE}s driven by Brownian
motions, L\'evy motions and \textsc{SDE}s driven by L\'evy motions.
For L\'evy motions, in particular, we consider symmetric $\a$-stable
L\'evy motions. We will see that the escape probability can be found by solving deterministic partial differential equations or partial differential-integral equations, with properly chosen boundary conditions.

\subsection{Boundary value problem  for escape probability of Brownian motions}

\begin{figure}[!htb]
 \centerline{\scalebox{1}{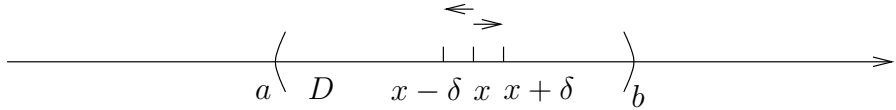}}
 \caption{A particle executing unbiased random walk in a bounded
interval}\label{fig1}
\end{figure}

Suppose that a particle executes an unbiased random walk on a
straight line. Let $D=(a,b)$. Figure \ref{fig1} shows the random
walk scenario. That is, a particle moves according to the following
rules (\cite{clls}):

  (i)  During the passage of a certain fixed time interval, a particle
  takes 1 step of a certain fixed length $\delta$ along
  the $x$ axis.

  (ii) It is equally probable that the step is to the right or to
  the left.

If the particle starting from $x\in D$ eventually escapes $D$ by
crossing the boundary $b$, then it must have moved to one of the two
points adjacent to $x$ first and then crossed the boundary. Thus \ce
p(x)=\frac{1}{2}[p(x-\delta)+p(x+\delta)], \de for $x\in D$. By
Taylor expansion on the right hand side to the second order, we have
$$
\frac{1}{2}p''(x)=0.
$$
The boundary conditions are
$$
\lim_{x\rightarrow b}p(x)=1, \quad\lim_{x\rightarrow a}p(x)=0,
$$
since the nearer the particle starts to $b$, the more likely
  it will first cross the boundary through $b$.

Note that the limit of the random walk is a standard Brownian motion $W_t$, that is,

(i) $W$ has independent increments;

(ii) for $0<s<t$, $W_t-W_s$ is a Gaussian random variable with mean zero
and variance $(t-s)$.\\
Thus, the escape probability $p(x)$ of a standard Brownian motion from $D$
through the boundary $b$ satisfies
\ce
\left\{\begin{array}{l}
\frac12 \Delta p(x)=0,\\
p(b)=1,\\
p(a)=0,
\end{array}
\right.
\de
where $\frac12 \Delta =\frac12 \partial_{xx}$ is the infinitesimal  generator for a scalar
standard Brownian motion $W_t$.

\subsection{Boundary value problem for escape probability of \textsc{SDE}s driven by Brownian motions}
\label{sde-bm}

Some results in this subsection can be found in \cite[Chapter 9]{ok}.

Let $\{W(t)\}_{t\geq 0}$ be an $m$-dimensional standard
$\cF_t$-adapted Brownian motion. Consider the following stochastic
differential equation (\textsc{SDE}) in $\mR^d$:
\be
X_{t}(x)=x+\int^{t}_{0}b(X_{s}(x))\,\dif
s+\int^{t}_{0}\sigma(X_{s}(x))\,\dif W_{s}. \label{eq}
\ee
We make the following assumptions about     the drift
$b:\mR^d\mapsto\mR^d$ and the diffusion coefficient $\sigma:\mR^d\mapsto\mR^d\times\mR^m$.

{\bf(H$^1_{b,\sigma}$)}
\ce
|b(x)-b(y)|\leq\lambda(|x-y|),\\
|\sigma(x)-\sigma(y)|\leq\gamma(|x-y|).
\de
Here $\lambda$ and $\gamma$ are increasing concave functions with the  properties
$\lambda(0)=\gamma(0)=0$, $\int_{0+}\frac{1}{\lambda(u)}\dif u=\int_{0+}\frac{1}{\gamma^2(u)}\dif u=\infty$.

\begin{figure}[!htb]
\centerline{\scalebox{1}{
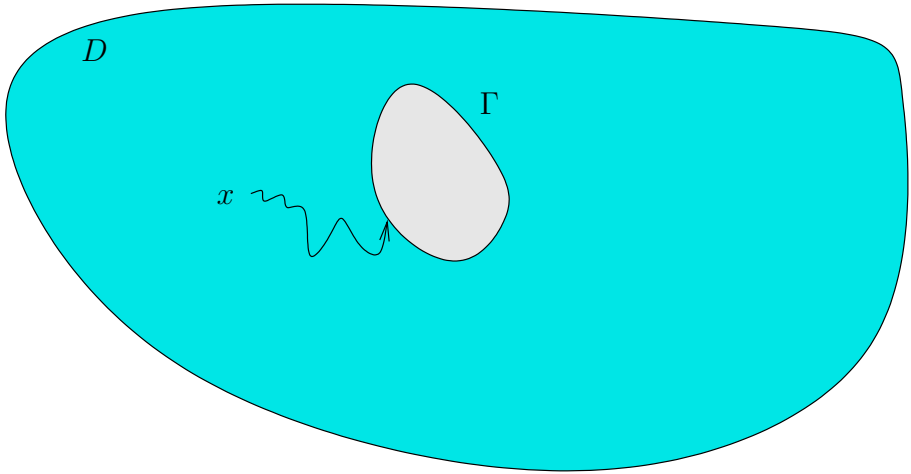}}
\caption{Escape probability for \textsc{SDE}s driven by Brownian motions: An annular open domain $D$ with a subset  $\Gamma$ of its boundary $\p D$} \label{fig2}
\end{figure}

Under {\bf(H$^1_{b,\sigma}$)}, it is well known that there exists a
unique strong solution to Eq.(\ref{eq})(\cite{tysw}). This solution is
denoted by $X_t(x)$.

We also make the following assumption.

{\bf(H$^2_{\sigma}$)}
There exists a $\xi>0$ such that for any $x,y\in D$
$$
\<y,\sigma\sigma^*(x)y\>\geq\xi|y|^2.
$$

This condition guarantees that the infinitesimal generator
\ce
\cL:=\sum\limits_{i=1}^db_i(x)\frac{\partial}{\partial x_i}+\sum\limits_{i,j=1}^da_{ij}(x)
\frac{\partial^2}{\partial x_i\partial x_j}
\de
for Eq.\eqref{eq} is uniformly elliptic in $D$, since then the eigenvalues of $\sigma\sigma^*$ are away
from $0$ in $D$. Here the matrix $[a_{ij}]:=\frac{1}{2}\sigma(x)\sigma^*(x)$.

Let $D$ be an open annular domain as in Figure \ref{fig2}. In one dimensional case, it is just an open interval. Let $\Gamma$ be its inner (or outer) boundary.  Taking
 \begin{eqnarray} \label{psi888}
\psi(x)=\left\{\begin{array}{l}
1, \quad x\in \Gamma,\\
0, \quad x\in \partial D\setminus\Gamma,
\end{array}
\right.
\end{eqnarray}
we have
\ce
\mE[\psi(X_{\tau_{\partial D}}(x))]&=&\int_{\{\omega: X_{\tau_{\partial D}}(x)\in\Gamma\}}\psi(X_{\tau_{\partial D}}(x))\dif \PX(\omega)\\
&&+\int_{\{\omega: X_{\tau_{\partial D}}(x)\in\partial D\setminus\Gamma\}}\psi(X_{\tau_{\partial D}}(x))\dif \PX(\omega)\\
&=& \PX\{\omega: X_{\tau_{\partial D}}(x)\in\Gamma\}\\
&=&p(x).
\de
This means that, for this specific $\psi$, $\mE[\psi(X_{\tau_{\partial D}}(x))] $ is the escape probability $p(x)$, which we are looking for.

We need to use \cite[Theorem 9.2.14]{ok} or \cite{ChenZQ} in order to see that the escape probability $p(x)$ is closely
related to a harmonic function with respect to $X_t$. This requires that the boundary data $\psi$ to be bounded
and continuous on $\p D$. For the domain $D$ taken as in Figure \ref{fig2}, with $\Gamma$ the inner
boundary (or outer) boundary, the above chosen $\psi$ in \eqref{psi888} is indeed bounded and continuous on $\p D$.
Thus, we have  the  following result by \cite[Theorem 9.2.14]{ok}.

\begin{theorem}
The escape probability $p(x)$ from an open annular domain $D$ to its inner (or outer)
boundary $\Gamma$, for the dynamical system driven by Brownian motions \eqref{eq}, is
the solution to the following Dirichlet boundary value problem
\ce
\left\{\begin{array}{l}
\cL p=0,\\
p|_{\Gamma}=1,\\
p|_{\partial D\setminus\Gamma}=0.
\end{array}
\right.
\de
\end{theorem}

\subsection{Boundary value problem for escape probability of symmetric $\alpha$-stable L\'evy motions}

\begin{figure}[!htb]
\centerline{\scalebox{1}{
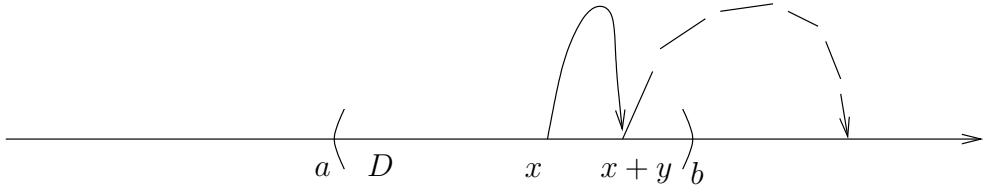}}
\caption{A particle executing L\'evy motion in   a bounded
interval} \label{fig3}
\end{figure}

Assume a particle is taking a one-dimensional L\'evy flight, where the distribution of
step sizes is a symmetric $\alpha$-stable distribution (Figure \ref{fig3}). Let $p(x)$ denote the escape probability of the particle starting
at $x$ in  $D=(a, b)$ and then first escapes $D$ over the right boundary $b$.
It could first move to somewhere inside $D$, say $x+y\in D$, and then
achieve its goal by jumping over the right boundary $b$ from the new starting point
$x+y$. More precisely,
\be\label{eq1}
p(x)=\int_{\mR\setminus\{0\}}P\{\mbox{the first step length is}~y\}p(x+y)\dif y.
\ee

According to \cite{Bergstr}, the symmetric $\alpha$-stable probability density function is the
following
\ce
f_{\a,0}(y)=\begin{cases}
-\frac{1}{\pi}\sum_{k=1}^\infty\frac{(-1)^k}{k!}\frac{\Gamma(\a
k+1)}{y|y|^{\a k}}\sin[k(\frac{\a\pi}{2}-\a \arg y)],\quad 0<\a<1,\\
\frac{1}{\pi}\sum_{k=0}^\infty(-1)^k\frac{\Gamma(\frac{k+1}{\a})}{k!\a}y^k\cos[k(\frac{\pi}{2})],
~\qquad\qquad\quad 1<\a<2,
\end{cases}
\de
where $\arg y=\pi$ when $y<0$.

For $0<\a<2$, the asymptotic expansion has also been given by
\cite{Bergstr} as follows
\ce
f_{\a,0}(y)&=&\begin{cases}
-\frac{1}{\pi}\sum_{k=1}^n\frac{(-1)^k}{k!}\frac{\Gamma(\a
k+1)}{y|y|^{\a k}}\sin\left[k(\frac{\a\pi}{2}-\a \arg y)\right]+o(|y|^{-\a(n+1)-1}),\quad |y|\rightarrow\infty,\\
\frac{1}{\pi}\sum_{k=0}^n(-1)^k\frac{\Gamma(\frac{k+1}{\a})}{k!\a}y^k\cos\left[k(\frac{\pi}{2})\right]
+o(|y|^{n+1}),\qquad\quad\qquad\quad\quad\quad |y|\rightarrow0,
\end{cases}\\
&=&\begin{cases} C_1(\a)/|y|^{1+\a}+o(|y|^{-2\a-1}),\quad |y|\rightarrow\infty,\\
C_2(\a)+o(|y|^2),\quad\qquad\quad\quad\quad |y|\rightarrow0,
\end{cases}
\de
where $C_1(\a)=\frac{1}{\pi}\sin(\frac{\pi\a}{2})\Gamma(1+\a)$ and
$C_2(\a)=\frac{1}{\pi}\frac{\Gamma(1/\a)}{\a}$. Take $N>0$ large enough and fix it. Thus,
\ce
0&=&\int_{\mR\setminus\{0\}}f_{\a,0}(y)[p(x+y)-p(x)]\dif y\\
&=&\int_{(-N,N)\setminus\{0\}}f_{\a,0}(y)[p(x+y)-p(x)]\dif y\\
&&+\int_{\mR\setminus(-N,N)}f_{\a,0}(y)[p(x+y)-p(x)]\dif y\\
&=:&I_1+I_2.
\de
For $I_1$, by self-affine property in \cite{Viswanathan}, we obtain
\be\label{eq2}
I_1&=&\int_{(-N,N)\setminus\{0\}}f_{\a,0}(N\cdot\frac{y}{N})[p(x+y)-p(x)]\dif y\no\\
&=&\int_{(-N,N)\setminus\{0\}}\frac{f_{\a,0}(N)}{(|y|/N)^{1+\a}}[p(x+y)-p(x)]\dif y\no\\
&=&\int_{(-N,N)\setminus\{0\}}\frac{C_1(\a)}{N^{1+\a}(|y|/N)^{1+\a}}[p(x+y)-p(x)]\dif y\no\\
&=&\int_{(-N,N)\setminus\{0\}}\frac{C_1(\a)}{|y|^{1+\a}}[p(x+y)-p(x)]\dif y.
\ee

For $I_2$, we calculate
\be\label{eq3}
I_2
&=&\int_{N}^{\infty}f_{\a,0}(y)[1-p(x)]dy+\int_{-\infty}^{-N}f_{\a,0}(y)[0-p(x)]\dif y\no\\
&=&\int_{N}^{\infty}\bigg[\frac{C_1(\a)}{y^{1+\a}}+o\bigg(\frac{1}{y^{1+2\a}}\bigg)\bigg][1-p(x)]\dif y\no\\
&&\qquad\qquad-\int_{-\infty}^{-N}\bigg[\frac{C_1(\a)}{(-y)^{1+\a}}+o\bigg(\frac{1}{(-y)^{1+2\a}}\bigg)\bigg]p(x)\dif y\no\\
&=&\int_{N}^{\infty}\frac{C_1(\a)}{y^{1+\a}}\dif y-\int_{\mR\setminus[-N,N]}\frac{C_{1}(\a)}{|y|^{1+\a}}p(x)\dif y\no\\
&=&\int_{\mR\setminus[-N,N]}\frac{C_1(\a)p(x+y)}{|y|^{1+\a}}\dif y-\int_{\mR\setminus[-N,N]}\frac{C_{1}(\a)}{|y|^{1+\a}}p(x)\dif y
\ee

Note that for $0<\a<1$, by the fact that the integral of an odd function on a symmetric
interval is zero, it holds that
\be
\int_{\mR\setminus\{0\}}p'(x)yI_{\{|y|\leq1\}}\frac{C_1(\a)}{|y|^{1+\a}}\dif y
=p'(x)C_1(\a)\int_{\{|y|\leq1\}\setminus\{0\}}\frac{y}{|y|^{1+\a}}\dif y
=0.
\label{intzero}
\ee
Thus, putting \eqref{eq2}, \eqref{eq3} and \eqref{intzero} together, we have for $0<\a<1$
\ce
\int_{\mR\setminus\{0\}}\left[p(x+y)-p(x)-p'(x)yI_{\{|y|\leq1\}}\right]\frac{C_1(\a)}{|y|^{1+\a}}\dif y=0.
\de
Moreover, $C_1(\a)=C_{1,\a}$.

For $\alpha\in[1,2)$, we only divide $I_1$ into two parts $I_{11}$ and $I_{12}$,
where
\ce
I_{11}&:=&\int_{\{|y|\leq\varrho\}\setminus\{0\}}f_{\a,0}(y)[p(x+y)-p(x)]\dif y,\\
I_{12}&:=&\int_{(-N,N)\setminus(-\varrho,\varrho)}f_{\a,0}(y)[p(x+y)-p(x)]\dif y,
\de
and $\varrho>0$ is a small enough constant.

For $I_{11}$, by Taylor expansion and self-affine property in [26],
we get
\ce
&&\int_{\{|y|\leq\varrho\}\setminus\{0\}}f_{\a,0}(y)[p(x+y)-p(x)]\dif y\\
&=&\int_{\{|y|\leq\varrho\}\setminus\{0\}}f_{\a,0}(\frac{1}{\varrho}\cdot\varrho y)p'(x)y\dif y\\
&=&\int_{\{|y|\leq\varrho\}\setminus\{0\}}\frac{f_{\a,0}(\frac{1}{\varrho})}
{(\varrho|y|)^{1+\alpha}}p'(x)y\dif y\\
&=&\int_{\{|y|\leq\varrho\}\setminus\{0\}}\frac{C_1(\a)}{(\frac{1}{\varrho})^{1+\alpha}\cdot(\varrho|y|)
^{1+\alpha}}p'(x)y\dif y\\
&=&\int_{\mR\setminus\{0\}}p'(x)yI_{\{|y|\leq\varrho\}}\frac{C_1(\a)}{|y|^{1+\alpha}}\dif y.
\de

For $I_{12}$, we apply the same technique as that in dealing with $I_1$ for $\alpha\in(0,1)$.

Next, by the similar calculation to that for $\alpha\in(0,1)$, we obtain for $\alpha\in[1,2)$
\ce
\int_{\mR\setminus\{0\}}\left[p(x+y)-p(x)-p'(x)yI_{\{|y|\leq1\}}\right]\frac{C_1(\a)\dif y}{|y|^{1+\a}}=0.
\de

Since the limit of the L\'evy flight is a symmetric $\a$-stable L\'evy motion $L_t^{\a}$, the
escape probability $p(x)$ of a symmetric $\a$-stable L\'evy motion, from $D$
to $[b,\infty)$ satisfies
\ce
\left\{\begin{array}{l}
-(-\Delta)^{\frac{\a}{2}}p(x)=0,\\
p(x)|_{[b,\infty)}=1,\\
p(x)|_{(-\infty,a]}=0.
\end{array}
\right.
\de
Note that $-(-\Delta)^{\frac{\a}{2}}$ is the infinitesimal generator for a scalar  symmetric $\a$-stable L\'evy motion $L_t^{\a}$.

\subsection{Boundary value problem for escape probability of \textsc{SDE}s driven by general L\'evy motions}
\label{sde-levy}

Let $L_t$ be a L\'evy process independent of $W_t$. Consider the following SDE in $\mR^d$
\be
X_{t}(x)=x+\int^{t}_{0}b(X_{s}(x))\,\dif s+\int^{t}_{0}\sigma(X_{s}(x))\,\dif W_{s}
+L_t.
\label{eqj1}
\ee
Assume that the  drift $b$ and the  diffusion $\sigma$ satisfy the following conditions:

{\bf (H$_b$)}
there exists a constant $C_b>0$ such that for $x,y\in\mR^{d}$
$$
|b(x)-b(y)|\leq C_b|x-y|\cdot
\log(|x-y|^{-1}+\me);
$$

{\bf (H$_\sigma$)}
there exists a constant $C_\sigma>0$ such that for $x,y\in\mR^{d}$
$$
|\sigma(x)-\sigma(y)|^2\leq C_\sigma|x-y|^{2}\cdot
\log(|x-y|^{-1}+\me).
$$

Under {\bf (H$_b$)} and {\bf (H$_\sigma$)}, it is well known that there exists a
unique strong solution to Eq.(\ref{eqj1})(see \cite{hqxz}). This solution will be
denoted by $X_t(x)$. Moreover, $X_t(x)$ is continuous in $x$.

\bl\label{sm}
The solution process $X_t(x)$ of the SDE \eqref{eqj1}  is a strong Markov process.
\el
\begin{proof}
Let $\eta$ be a $(\cF_t)_{t\geq0}$-stopping time. Set
\ce
\cG_t:=\sigma\{W_{\eta+t}-W_\eta, L_{\eta+t}-L_\eta\}\cup\cN, \quad t\geq0,
\de
where $\cN$ is of all $P$-zero sets. That is, $\cG_t$ is a completed $\sigma$-algebra
generated by $W_{\eta+t}-W_\eta$ and $L_{\eta+t}-L_\eta$. Besides, $\cG_t$ is independent
of $\cF_t$. Let $X(x,\eta,\eta+t)$ denote the unique solution of the the following SDE
\be
X(x,\eta,\eta+t)=x+\int^{\eta+t}_{\eta}b(X(x,\eta,s))\,\dif s+\int^{\eta+t}_{\eta}\sigma(X(x,\eta,s))\,\dif W_{s}
+L_{\eta+t}-L_{\eta}.
\label{eqj2}
\ee
Moreover, $X(x,\eta,\eta+t)$ is $\cG_t$-measurable and $X(x,0,t)=X_t(x)$. By the
uniqueness of the solution to (\ref{eqj2}), we have
\ce
X(x,0,\eta+t)=X(X(x,0,\eta),\eta,\eta+t), \quad a.s..
\de

For any bounded measurable function $g$,
\be
\mE[g(X_{\eta+t}(x))|\cF_\eta]&=&\mE[g(X(x,0,\eta+t))|\cF_\eta]\no\\
&=&\mE[g(X(X(x,0,\eta),\eta,\eta+t))|\cF_\eta]\no\\
&=&\mE[g(X(y,\eta,\eta+t))]|_{y=X(x,0,\eta)}\no\\
&=&\mE[g(X(y,0,t))]|_{y=X(x,0,\eta)}.
\label{smp}
\ee
Here the last equality holds because the distribution of $X(y,\eta,\eta+t)$ is the
same to that of $X(y,0,t)$. The proof is completed since (\ref{smp}) implies that
\ce
\mE[g(X_{\eta+t}(x))|\cF_\eta]=\mE[g(X_{\eta+t}(x))|X_{\eta}(x)].
\de
\end{proof}

\begin{figure}[!htb]
\centerline{\scalebox{1}{
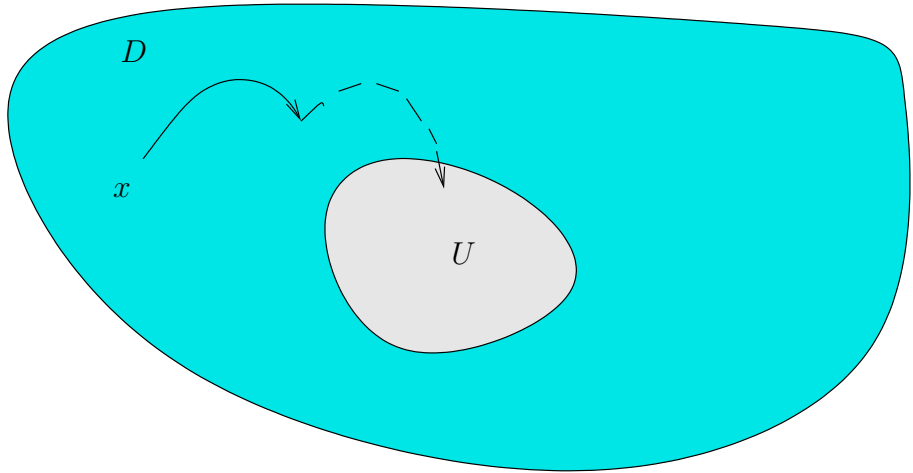}}
\caption{Escape probability for \textsc{SDE}s driven by L\'evy
motions: an open annular domain $D$, with its inner part $U$
(which is in $D^c$) as a target domain } \label{fig4}
\end{figure}

\begin{figure}[!htb]
\centerline{\scalebox{1}{
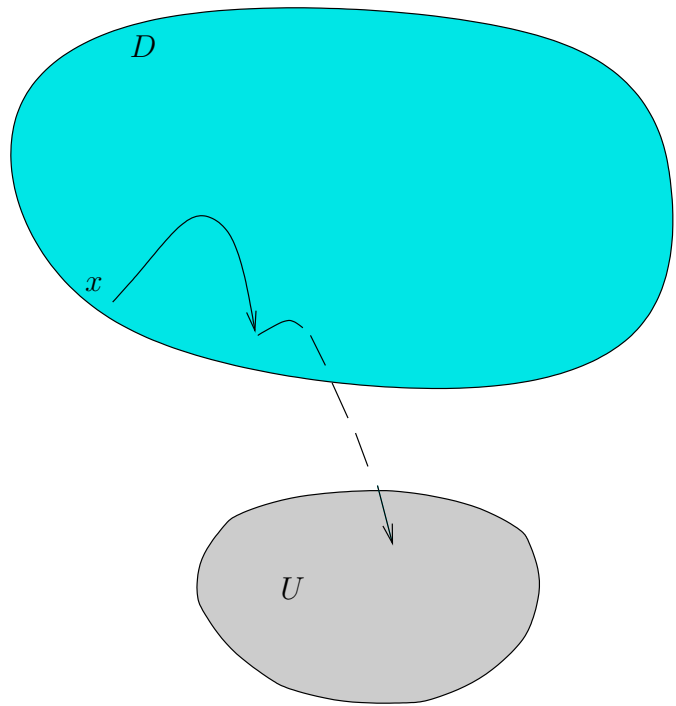}}
\caption{Escape probability for \textsc{SDE}s driven by L\'evy
motions: a general open domain $D$, with a target domain $U$  in $D^c$}
\label{fig5}
\end{figure}

Because $L_t$ has c\`adl\`ag and quasi-left-continuous paths(\cite{s}), $X_t(x)$ also has c\`adl\`ag and quasi-left-continuous paths. Thus
by Lemma \ref{sm} and Definition \ref{hunt process} above, we see that $X_t(x)$
is a Hunt process. Let $D$ be a relatively compact and regular open domain (Figure \ref{fig4} or Figure \ref{fig5}). Theorem \ref{soldir}
implies that $\mE[\varphi(X_{\tau_{D^c}}(x))]$ is the unique solution to
the Balayage-Dirichlet problem $(D, \varphi)$, under the condition that $\varphi$ is nonnegative and bounded on $D^c$.
Set
 \ce
\varphi(x)=\left\{\begin{array}{l}
1, \quad x\in U,\\
0, \quad x\in  D^c \setminus U.
\end{array}
\right.
\de
Then $\varphi$ is nonnegative and bounded on $D^c$. We observe that
\ce
\mE[\varphi(X_{\tau_{D^c}}(x))]&=&\int_{\{\omega: X_{\tau_{D^c}}(x)\in U\}}\varphi(X_{\tau_{D^c}}(x))\dif \mP(\omega)\\
&&+\int_{\{\omega: X_{\tau_{D^c}}(x)\in D^c\setminus U\}}\varphi(X_{\tau_{D^c}}(x))\dif \mP(\omega)\\
&=&\mP\{\omega: X_{\tau_{D^c}}(x)\in U\}\\
&=&p(x).
\de
This means that, for this specific $\varphi$,  $\mE[\varphi(X_{\tau_{D^c}}(x))]$ is the escape probability $p(x)$ that we are looking for.
By the definition of the characteristic operator, $p\in I_D$ and by Theorem \ref{soldir}, $Ap(x)=0$. Thus we obtain the following theorem.

\bt
Let $D$ be a relatively compact
and regular open domain, and let $U$ be a set in $D^c$. Then the escape probability $p(x)$, for the dynamical system driven by L\'evy motions \eqref{eqj1}, from $D$ to $U$,
is the solution of the following Balayage-Dirichlet problem
\ce
\left\{\begin{array}{l}
A p=0,\\
p|_{U}=1,\\
p|_{D^c\setminus U}=0,
\end{array}
\right.
\de
where $A$ is the characteristic operator for this system.
\et

\begin{remark}
Unlike the SDEs driven by Brownian motions, a typical open domain $D$ here could be a quite general open   domain (Figure \ref{fig5}), as well as   an annular domain (Figure \ref{fig4}).
This is due to the jumping properties of the solution paths. It is also due to the fact that, in Theorem \ref{soldir}, the   function $f$ is only required to be continuous on the boundary $\p D$ (not on the domain $D^c$).
\end{remark}

\bigskip
\medskip

Finally we consider the representation of the characteristic
operator $A$,  for   a SDE driven by a symmetric $\alpha$-stable
L\'evy process $L_t^{\alpha}$, with $\alpha \in (0, 2)$:
\be
X_{t}(x)=x+\int^{t}_{0}b(X_{s}(x))\,\dif
s+\int^{t}_{0}\sigma(X_{s}(x))\,\dif W_{s} +L_t^{\alpha}.
\label{sdelevy999}
\ee

 Let us first consider the case of $1\leq\a<2$. For $f\in\cC_0^2(\mR^{d})$, applying the
It\^o formula to $f(X_{\tau_{\varepsilon}}(x))$,  we obtain
\ce
f(X_{\tau_{\varepsilon}}(x))-f(x)&=&\int_0^{\tau_{\varepsilon}}\<\partial_yf(X_s),b(X_s)\>\dif s
+\int_0^{\tau_{\varepsilon}}\<\partial_yf(X_s),\sigma(X_s)\dif W_s\>\\
&&+\int_0^{\tau_{\varepsilon}}\int_{|u|\leq1}\left(f(X_s+u)-f(X_s)\right)\tilde{N}_{p}(\dif s, \dif u)\\
&&+\int_0^{\tau_{\varepsilon}}\int_{|u|>1}\left(f(X_s+u)-f(X_s)\right)N_{p}(\dif s, \dif u)\\
&&+\frac{1}{2}\int_0^{\tau_{\varepsilon}}\left(\frac{\partial^2}{\partial y_i\partial y_j}f(X_s)\right) \sigma_{ik}(X_s)\sigma_{kj}(X_s) \dif s\\
&&+\int_0^{\tau_{\varepsilon}}\int_{|u|\leq1}\big(f(X_s+u)-f(X_s)\\
&&\qquad-\<\partial_yf(X_s),u\>\big)
\frac{C_{d,\alpha}}{|u|^{d+\alpha}}\dif u\dif s.
\de
Here and hereafter, we use the convention that repeated indices imply summation from $1$ to $d$.
Taking expectation on both sides, we get
\ce
&&\mE[f(X_{\tau_{\varepsilon}}(x))]-f(x)\\
&=&\mE\int_0^{\tau_{\varepsilon}}\<\partial_yf(X_s),b(X_s)\>\dif s
+\frac{1}{2}\mE\int_0^{\tau_{\varepsilon}}\left(\frac{\partial^2}{\partial y_i\partial y_j}f(X_s)\right)\sigma_{ik}(X_s)\sigma_{kj}(X_s)\dif s\\
&&+\mE\int_0^{\tau_{\varepsilon}}\int_{\mR^d\setminus\{0\}}\big(f(X_s+u)-f(X_s)-\<\partial_yf(X_s),u\>\big)
\frac{C_{d,\alpha}}{|u|^{d+\alpha}}\dif u\dif s.
\de

The infinitesimal generator $\cL$ of Eq.(\ref{eqj1}) is as follows \cite{da}:
\ce
(\cL f)(x)&=&\<\partial_xf(x),b(x)\>+\frac{1}{2}\left(\frac{\partial^2}{\partial x_i\partial x_j}
f(x)\right) \sigma_{ik}(x) \sigma_{kj}(x)\\
&&+\int_{\mR^d\setminus\{0\}}\left(f(x+u)-f(x)
-\<\partial_xf(x),u\>\right)
\frac{C_{d,\alpha}}{|u|^{d+\alpha}}\dif u.
\de
So,
\ce
\left|\frac{\mE[f(X_{\tau_{\varepsilon}}(x))]-f(x)}{\mE[\tau_{\varepsilon}]}-(\cL f)(x)\right|
&=&\left|\frac{\mE\int_0^{\tau_{\varepsilon}}(\cL f)(X_s)\dif s}{\mE[\tau_{\varepsilon}]}
-\frac{\mE\int_0^{\tau_{\varepsilon}}(\cL f)(x)\dif s}{\mE[\tau_{\varepsilon}]}\right|\\
&\leq&\frac{\mE\int_0^{\tau_{\varepsilon}}|(\cL f)(X_s)-(\cL f)(x)|\dif s}{\mE[\tau_{\varepsilon}]}\\
&\leq&\sup\limits_{|y-x|<\e}|(\cL f)(y)-(\cL f)(x)|.
\de
Because $(\cL f)(x)$ is continuous in $x$,
\ce
A f(x)=\lim\limits_{\e \downarrow 0}\frac{\mE[f(X_{\tau_{\varepsilon}}(x))]
-f(x)}{\mE[\tau_{\varepsilon}]}=(\cL f)(x).
\de

Similarly, we also have $A=\cL$ for $0<\a<1$.

\medskip

\br
The above deduction tells us $Af=\cL f$ for $f\in\cC_0^2(\mR^{d})$. If the considered driving process
is not   a symmetric $\alpha$-stable L\'evy motion, the domain of $\cL$ is unclear and thus     $A=\cL$ may not be true. The corresponding
escape probability $p(x)$ is the solution of the following Balayage-Dirichlet problem (in terms of operator $\cL$, instead of $A$):
$$
\left\{\begin{array}{l}
\mathcal {L}p=0,\\
p|_{U}=1,\\
p|_{D^c\setminus U}=0.
\end{array}
\right.
$$
\er

\section{Examples}\label{example}

In this section we consider a few examples.

\begin{example} \label{example4.1}
 In 1-dimensional case, take $D=(-r,r)$ and $\Gamma=\{r\}$. For each $x\in D$, the escape
probability $p(x)$ of $X_t=x+W_t$ from $D$ to $\Gamma$ satisfies the
following differential equation
\ce\left\{\begin{array}{l}
\frac{1}{2}p''(x)=0, \quad x\in (-r,r),\\
p(r)=1,\\
p(-r)=0.
\end{array}
\right.
\de
We obtain that $p(x)=\frac{x+r}{2r}$ for $x\in[-r,r]$. It is a straight line (See Figure \ref{fig6}).

\begin{figure}
{\scalebox{0.7}{\includegraphics*[25,30][550,580]{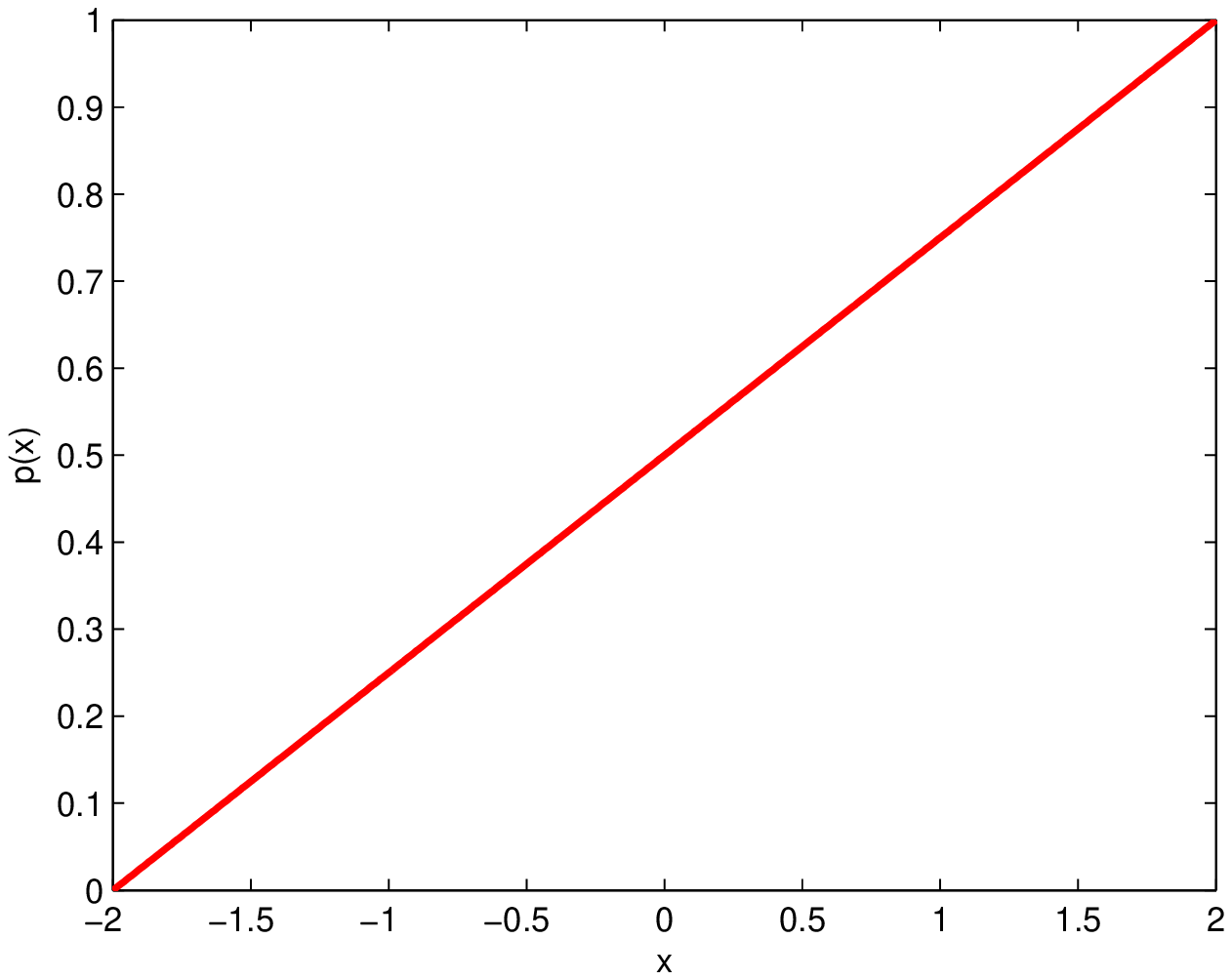}}}\\

\vspace*{-55mm}

\caption{Escape probability for 1-dimensional Brownian motion in Example \ref{example4.1}, $r=2$} \label{fig6}
\end{figure}

In 2-dimensional case, take $D=\{x\in\mR^2; r<|x|<R\}$ and $\Gamma=\{x\in\mR^2; |x|=r\}$. For every $x\in D$, the escape
probability $p(x)$ of $X_t=x+W_t$ from $D$ to $\Gamma$ satisfies the
following elliptic partial differential equation
\ce\left\{\begin{array}{l}
\frac{1}{2}\Delta p(x)=0, \quad x\in D,\\
p(x)|_{|x|=r}=1,\\
p(x)|_{|x|=R}=0.
\end{array}
\right.
\de
By solving this equation, we obtain that
$p(x)=\frac{\log R-\log|x|}{\log R-\log r}$. It is plotted in Figure \ref{fig7}.

\begin{figure}
{\scalebox{0.7}{\includegraphics*[25,30][550,580]{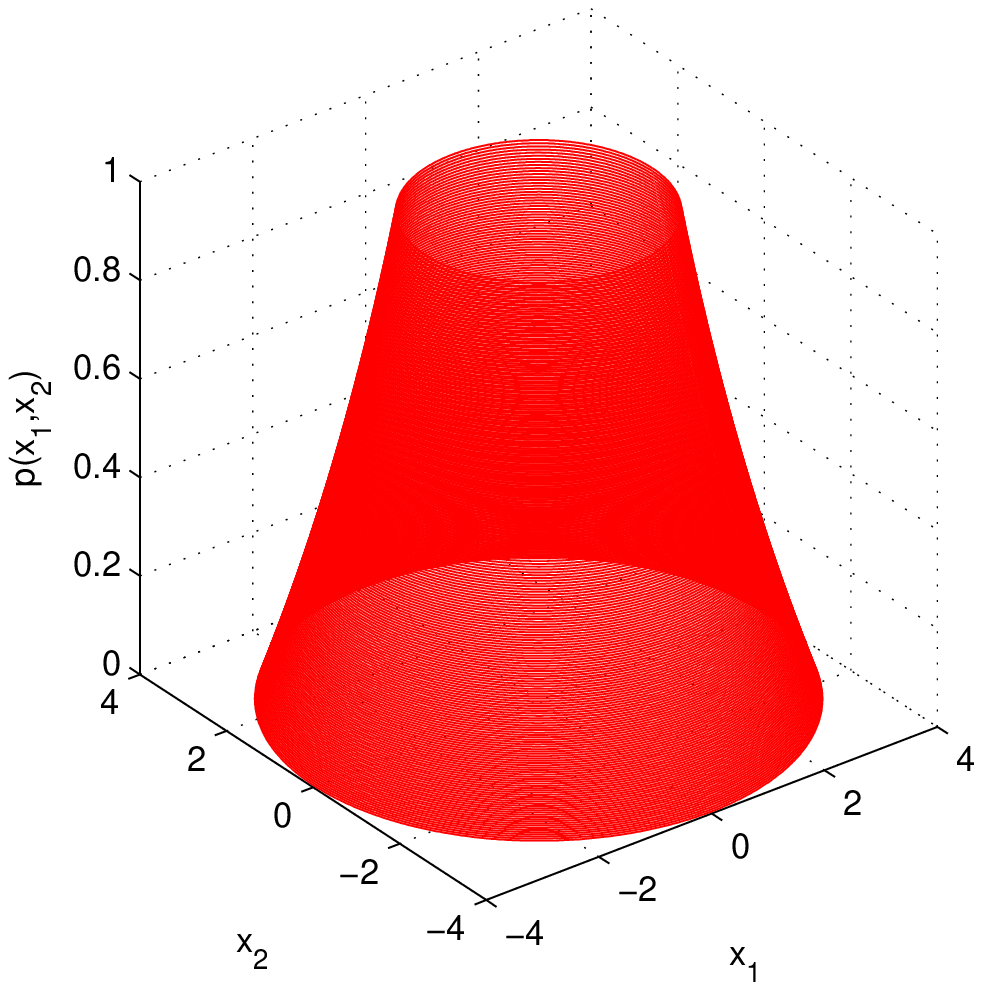}}}\\

\vspace*{-55mm}

\caption{Escape probability for 2-dimensional Brownian motion in Example \ref{example4.1}: $r=2, R=4$} \label{fig7}
\end{figure}
\end{example}

\begin{example}\label{example4.2}
 Consider the following \textsc{SDE} driven by Brownian motions:
 \ce
\dif X_t=b(X_t)\dif t+\sigma(X_t)\dif W_t,
\de
where $b$ and (non-zero)
$\sigma$ are real functions. When $b$ and $\sigma$ satisfy {\bf(H$^1_{b,\sigma}$)},
the equation has a unique solution which is denoted as $X_t$. We take $D=(-r,r)$
and $\Gamma=\{r\}$. For each $x\in D$, under the condition {\bf (H$^2_\sigma$)},
the escape probability $p(x)$ satisfies
$$
\left\{\begin{array}{l}
\frac{1}{2}\sigma^2(x)p^{''}(x)+b(x)p^{'}(x)=0, \quad x\in(-r,r),\\
p(r)=1,\\
p(-r)=0.
\end{array}
\right.
$$
The solution is
$$
p(x)=\frac{\int_{-r}^xe^{-2\int_{-r}^y\frac{b(z)}{\sigma^2(z)}\dif z}\dif y}
{\int_{-r}^re^{-2\int_{-r}^y\frac{b(z)}{\sigma^2(z)}\dif z}\dif y}
$$
for $x\in[-r,r]$.  See Figure \ref{fig8}.

\begin{figure}
{\scalebox{0.7}{\includegraphics*[25,30][550,580]{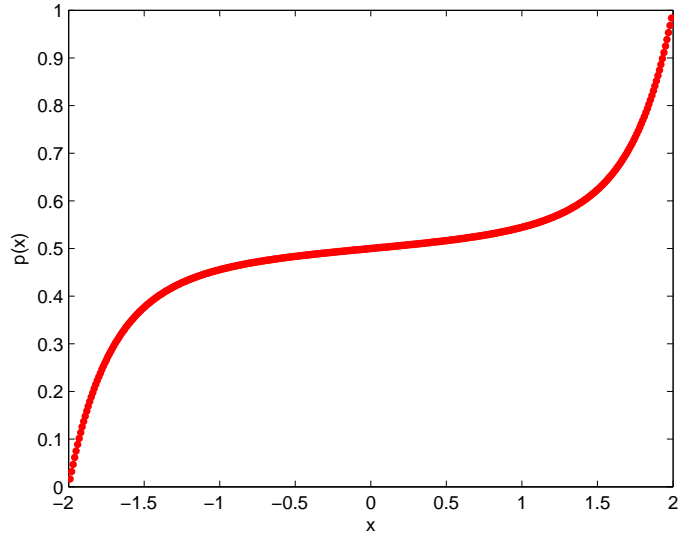}}}\\

\vspace*{-55mm}

\caption{Escape probability in Example \ref{example4.2}: $b(x)=-x, \sigma(x)=1, r=2$} \label{fig8}
\end{figure}
\end{example}

\begin{example} \label{example4.3}
 In 1-dimensional case, take $D=(-r,r)$ and $U=[r,\infty)$. For each $x\in D$
and a symmetric $\a$-stable L\'evy process $L_t^\alpha$, the escape
probability $p(x)$ of $X_t=x+L_t^\alpha$ from $D$ to $U$ satisfies the
following differential-integral equation
\ce\left\{\begin{array}{l}
-(-\Delta)^{\frac{\alpha}{2}}p(x)=0, \quad x\in(-r,r),\\
p(x)|_{[r,\infty)}=1,\\
p(x)|_{(-\infty,-r]}=0.
\end{array}
\right.
\de
It is difficult to deal with   this equation because of the fractional Laplacian operator. But
we can solve it via Poisson kernel. From \cite{ja}, for $x\in(-r,r)$,
\ce
p(x)=\frac{\sin\frac{\pi\a}{2}}{\pi}\int_r^\infty\frac{(r^2-x^2)^{\a/2}}{(y^2-r^2)^{\a/2}}\frac{1}{(y-x)}\dif y.
\de

Obviously, $p(-r)=0$. To justify $p(r)=1$, we apply the substitution $y=(r^2-xv)(x-v)^{-1}$ to obtain
\ce
p(r)&=&\frac{\sin\frac{\pi\a}{2}}{\pi}\int_{-r}^r(r-v)^{\a-1}(r^2-v^2)^{-\frac{\a}{2}}\dif v\\
&=&\frac{\sin\frac{\pi\a}{2}}{\pi}\int_0^1(1-v)^{\frac{\a}{2}-1}v^{1-\frac{\a}{2}-1}\dif v\\
&=&\frac{\sin\frac{\pi\a}{2}}{\pi}B(\frac{\a}{2},1-\frac{\a}{2})\\
&=&\frac{\sin\frac{\pi\a}{2}}{\pi}\Gamma(\frac{\a}{2})\Gamma(1-\frac{\a}{2})\\
&=&1,
\de
where the Beta and Gamma functions and their properties are used in the last two steps. The escape probability $p(x)$
is plotted in Figure \ref{fig9} for various $\alpha$ values.

\begin{figure}
\begin{center}
\includegraphics{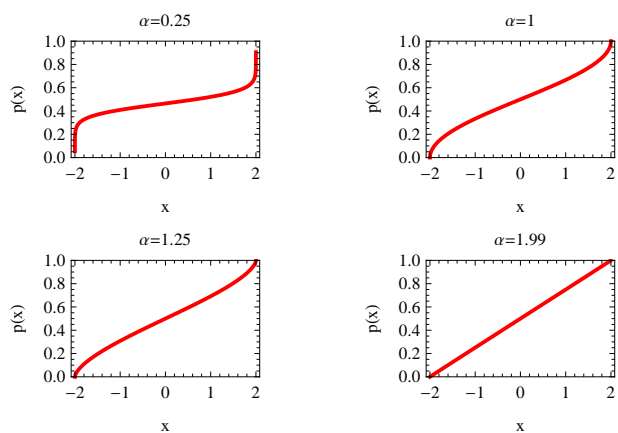}
\end{center}
\caption{Escape probability in Example \ref{example4.3}: $ r=2$} \label{fig9}
\end{figure}

\end{example}

\bigskip

\emph{Acknowledgement}.  We have benefited from our previous collaboration with Ting Gao, Xiaofan Li and Renming Song \cite{GaoDuan}. We   thank Ming Liao, Renming Song and Zhen-Qing Chen for helpful discussions. This work was done while Huijie Qiao was visiting the Institute for Pure and Applied Mathematics (IPAM), Los Angeles.

\end{document}